\documentclass[a4paper, 12pt]{amsart}  
\usepackage{times}
\usepackage{amssymb}
\usepackage[english]{babel}
\usepackage[autostyle]{csquotes}
\usepackage{mathrsfs}
\usepackage{amsmath,amssymb,amsfonts,enumerate,amsthm,graphicx,color}
\usepackage[all]{xy}
\theoremstyle{plain}
  \newtheorem{thm}{Theorem}[section]

\theoremstyle{definition}
  \newtheorem{dfn}[thm]{Definition}

\theoremstyle{remark}
  
\numberwithin{equation}{section}


\begin{document}
\title{ LOW CODIMENSION STRATA OF THE SINGULAR LOCUS OF  MODULI  OF LEVEL CURVES }
\author{ sepideh tashvighi}
\maketitle

\textbf{Abstract.}
We further analyse the moduli space of stable curves with level structure provided by Chiodo and Farkas in \cite{AA}. Their result builds upon Harris and Mumford analysis of the locus of singularities of the moduli space of curves and shows in particular that for levels 2, 3, 4, and 6 the locus of noncanonical singularities is completely analogous to the locus described by Harris and Mumford, it has codimension 2 and arises from the involution of elliptic tails carrying a trivial level structure. For the remaining levels (5, 7, and beyond), the picture also involves components of higher codimension. We show that there exists a component of codimension 3 for levels $\ell=5$ and $\ell\geqslant 7$ with the only exception of level 12. We also show  that there exists a component of codimension 4 for $\ell=12$.\\


\section{Introduction}
Let ${\mathcal M_{g}} $ be the moduli space of smooth curves of genus $g$. We denote by $\overline{\mathcal M}_{g}$ the compactification of $\mathcal M_{g}$. Its objects are nodal curves  satisfying Deligne-Mumford stability condition.
 Recall that  to a stable curve  one can attach a graph, the  so-called  stable graph.
 A graph is called stable if for all of its vertices $v$, the
following inequality  holds
 $2g_{v}-2+\ell_{v}>0 $, where $\ell_{v} $ is the valence, i.e the number of entering
edges and $g_v$ is the geometric genus of  the normalization of  the irreducible
component corresponding to $v$.\\
 Let $\mathcal{R}_{g,\ell}$ be the  moduli space  of  level structures of all triples
$(C,L,\phi)$, where $C$ is a smooth curve of genus $g$ equipped with a line
bundle $L$ and a trivialization morphism $\phi \colon L^{\otimes
\ell}\stackrel{\sim} \longrightarrow \mathcal{O}$. We refer to these as  level-$\ell$ curves.  We consider  $\overline{\mathcal R}_{g,\ell}$, a compactification of the moduli space $\mathcal R_{g,\ell}$. Similarly to the moduli space of stable
curves $\overline {\mathcal M}_{g} $ (see \cite{GG}), the locus $\overline{\mathcal
{R}}_{g,\ell}\setminus \mathcal R_{g,\ell} $ can be described by the
dual graph of nodal curves. We determine its vertices, $V$  by connected components and its edges, $E$  by  nodes of nodal curves. The dual graph $\Gamma$ is  a stable graph decorated by
${M=\lbrace m_{e}\rbrace }_{e \in{{E}}}$ ($ \mathbb{Z}_{\ell}$-valued 1-cochain from the set of branches of each node of $\mathsf{C}$)  lying in the kernel of the
 cochain homomorphism $\partial \colon
C^{1}(\Gamma,\mathbb{Z}_{\ell}){\longrightarrow}C^{0}(\Gamma,\mathbb{Z}_{\ell})$. This means that $M$ adds up to zero modulo $\ell$ along each circuit of $\Gamma$.\\
 Recall that $\overline{\mathcal M}_{g} $ is locally isomorphic to ${\mathrm{Def}(\mathsf{C},\mathsf{L},\phi)}\slash{\mathrm {Aut}(\mathsf{C})}$ where $\mathrm{Def}(\mathsf{C},\mathsf{L}, \phi)$ is the deformation space of the stack-theoretic curve. Note that a level structure $(\mathsf{C},\mathsf{L},\phi)$ is smooth if and only if each element of $\mathrm{Aut}(\mathsf{C})$ operates on $\mathrm{Def}(\mathsf{C},\mathsf{L},\phi)$ as a product of quasireflections (i.e., an automorphism whose fixed locus is hyperplane).\\
We consider the automorphisms  which are given by twisting each
node. If we  specify a coefficient in $\mathbb{Z}_{\ell}$ for every edge, then
every choice of coefficient determines an automorphism  $\mathsf{a}$ of the stack-theoretic curve $\mathsf{C}$. Note that the automorphisms
preserve the line bundle if and only if the action of the automorphism $\mathsf{a}$ on $M$ lies in  the image of $\delta \colon C^{0}(\Gamma,\mathbb{Z}_{\ell}){\longrightarrow}C^{1}(\Gamma,\mathbb{Z}_{\ell})$ (see \cite{AA}).\\
Then we consider the existence of noncanonical singularities of
${\mathrm{Def}(\mathsf{C},\mathsf{L},\phi)}\slash{\mathrm {Aut}(\mathsf{C})}$. It is  known that  the noncanonical singularities occur if and only if there exists a  junior automorphism (i.e., an automorphism \enquote{aged} less than one) on $\mathrm{Def}(\mathsf{C},\mathsf{L},\phi)$ over the group of automorphism mod out by quasireflections \cite{AA}.
Using this fact, we prove that there exists a codimension 3 locus of noncanonical singularities for all levels  $\ell=5$ and $\ell \geqslant 7$ with the only exception of level 12. We also show that there exists a codimension 4 locus of noncanonical singularities for $\ell=12$.

\section{Preliminaries}
 Assume that
$k$ is an algebraically closed field.
Assume $k$ is $\mathbb{C}$, the field of complex number.
 Let us recall some useful
definitions.

\begin{dfn} Let  $V$ be the set of connected components of normalization of stack-theoretic curve, $ \mathsf{C} $.
Let  $E$ be the set of nodes of $\mathsf{C}$. The pair $ \left( V,E\right) $
is called the dual graph of normalization of $\mathsf{C}$, i.e. $\mathsf{C}^{\nu}$ (see
 for instance \cite{CC}).
\end{dfn}
Now for every nodal curve we are able to draw its dual graph. We define $C^{0}(\Gamma,\mathbb{Z}_{\ell})=\lbrace a:V \rightarrow \mathbb{Z}
\rbrace=\bigoplus_{\upsilon \in V}\mathbb{Z}$ as the set of
$\mathbb{Z}$-valued functions on  $V$. We define $C^{1}(\Gamma,\mathbb{Z}_{\ell})=\lbrace b:{\mathbb{E}} \rightarrow
\mathbb{Z}  \mid b(\overline{e})=-b(e)\rbrace$ (i.e. ${\mathbb{E}}$ is the set
of branches of each node of $\mathsf{C}$)  as the set of
antisymmetric $\mathbb{Z}_{\ell}$-valued functions on $\mathbb{E}$, where
$\overline{e}$ and $e$ are oriented edges with opposite
orientations.

Let us recall that $\delta \colon
C^{0}(\Gamma,\mathbb{Z}_{\ell})\rightarrow
C^{1}(\Gamma,\mathbb{Z}_{\ell})$  is defined by sending $a$ to
$\delta a$, with $\delta a(e)=a(e_{+})-a(e_{-})$ and
  the map $\partial \colon C^{1}(\Gamma,\mathbb{Z}_{\ell}){\longrightarrow}C^{0}(\Gamma,\mathbb{Z}_{\ell})$ is
   defined by sending $b$ to $\partial b$, with $\partial b(v)={\Sigma}_{e \in {E}} b(e)$  \cite{AA}.

Let  $\mathsf{a}$ be an automorphism of stack-theoretic curve $\mathsf{C}$  and $M$ be the multiplicity cochain, then we can define $\mathsf{a}$ as a multiple of  $\gcd(M,\ell)$, where $\mathsf{a}$ is defined as  follows:
$$\mathsf{a} \odot M=\dfrac{{a}M}{\gcd(M,\ell)},$$
(see \cite[pages 35 and 36]{AA}).\\
We recall that, $\mathrm{age}(\mathsf{a})=\Sigma_{e \in
E} \left\lbrace \dfrac{\mathsf{a}(e)}{ \ell}  \right\rbrace $.\\
 \begin{dfn} Let $\mathrm{Aut}_{C}(\mathsf{C})$ be the group of automorphisms of $\mathsf{C}$.
 We say that an automorphism  $\mathsf{a}$ in  $\mathrm{Aut}_{C}(\mathsf{C})$, which  operates nontrivially on
 the curve, is junior on  $\mathsf{C}$ if $0<\mathrm{age}(\mathsf{a})<1 $,
 (see \cite[Definition 2.35]{AA}).
\end{dfn}

According to the    Reid--Shepherd-Barron--Tai criterion, the scheme theoretic-quotient $V/G$, where $G$ is a finite group operates on $V$ without quasireflection, has a noncanonical singularity at the origin if and only if the image of age$_{V}$  intersects $\left] 0,1 \right[ $, see \cite{EE, DD, FF}.
 As above, the point is that, the quotient
${\mathbb{C}^{3g-3}}\slash{\mathrm {Aut}(\mathsf{C})}$ has a
noncanonical singularities  if and only if there exists an element
$\mathsf{a} \in \mathrm{Aut}_{C}(\mathsf{C})$ which is junior on
$\mathsf{ C}$.
\begin{thm}(\cite{AA}) There exist a  junior
ghost $\mathsf{a}$ if the following conditions are satisfied:
\begin{itemize}
\item[(i)] $\mathrm{age}(\mathsf{a})<1$
(i.e., $\mathsf{a}$ is junior);
\item[(ii)] $M={\Sigma}_{i \in I} K_{i}$, where I is a finite set of
circuits (i.e., $M \in \mathrm{Ker}\ \partial$);
\item[(iii)] $\mathsf{a}\odot M(K)\equiv 0$ for any circuit K (i.e.
$\mathsf{a}\odot M \in \mathrm{Im} \ \delta$).
\end{itemize}
\end{thm}
\section{Classifying the noncanonical singularities}
In this section our aim is to analyse the existence of strata in the locus of noncanonical singularities of codimension 3. First of all, we show that there is a
codimension 3 locus of noncanonical singularities  for all levels 5,7, and higher except $\ell=12$. This is done in three steps.

Step1. Let  $\ell$ be a prime number bigger than 3. Consider a level curve whose dual graph has multiplicity $M$. Let  $m _1= m _3=n$ and
$m _2=2n$, where $n$
  can be chosen in $\mathbb{Z}/\ell$ (see Figure 1). Let $a_1=a_3=1$ and $a_2=\dfrac{\ell-1}{2}$ in $\mathsf{a}$.
  Hence the sum of the values of $a \odot M$ along every circuit is zero modulo $\ell$. Indeed $(a_{1}\odot m _{1})+(a_{2}\odot m _{2})$
  and $(a_{1}\odot m _{1})-(a_{3}\odot m _{3})$ add up to zero modulo $\ell$.\\
  $$
\begin{matrix}
\xymatrix{
\bullet \ar@/^2pc/[rr]^{m_1} \ar@/
_1.5pc/[rr]^{m_3}
& & \bullet \ar[ll]_{m_2}
}
\end{matrix}
$$
\begin{center}
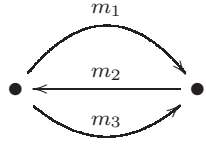
Figure 1. Graph $M$ with two vertices and three edges.
\end{center}

  On the other hand, since $\ell>3$, $\mathrm{age}(\mathsf{a})=\dfrac{a_{1}}{\ell}+\dfrac{a_{2}}{\ell}+\dfrac{a_{3}}{\ell}=\dfrac{\ell+3}{2\ell}<1$.

Notice that if we find a junior $\mathsf{a}$ our problem for a given $\ell$, then every
multiple can also be solved similarly. In fact, if for a given
$\ell$ we find $m_{1}, m_{2}$ and $m_{3}$ in $M$, as well as $a_{1},
a_{2}$ and $a_{3}$ in $\mathsf{a}$ which satisfy in our conditions,
then we can set $\ell^{'}=k\ell$, $m^{'}_{1}=km_{1},
m^{'}_{2}=km_{2}$ and $m^{'}_{3}=km_{3}$ in $M$, also
$a^{'}_{1}=ka_{1}, a^{'}_{2}=k a_{2}$ and $a^{'}_{3}=ka_{3}$ in
$\mathsf{a}$ for every integer $k$. Hence the sum along every circuit
is zero modulo $\ell^{'}$ and $(a^{'}_{1}\odot m^{'}
_{1})+(a^{'}_{2}\odot m^{'} _{2})$ , $(a^{'}_{1}\odot m^{'}
_{1})-(a^{'}_{3}\odot m^{'}_{3})$ add up to zero modulo $\ell^{'}$.

According to the above argument, there exists a codimension 3 locus of
noncanonical singularities for prime numbers $\ell>3$. This settles the cases $\ell\neq 2^{a}3^{b}$. We now focus on $\ell= 2^{a}3^{b}$.

Step 2: if $\ell=8$, take $m_{1}=1, m_{2}=3, m_{3}=2$ and
$a_{1}=a_{2}=a_{3}=2$. Then, there exists a junior ghost
automorphism.

If $\ell=9$, take $m_{1}=1, m_{2}=2, m_{3}=1$ and
$a_{1}=a_{3}=1$ and $a_{2}=4$. Then, there exists a junior ghost
automorphism.

Step 3. \cite{AA} show that there is no junior ghost $\mathsf{a}$ for any stable graph
for which $M \in \mathrm{Ker}\ \partial$ and $\mathsf{a} \odot M  \in
\mathrm{Im}\ \delta$ for $\ell=2, 3, 4, 6$. We also show that for $\ell=12$ if there are only 3 edges. To do this we provide a table.
In the first row, we write the  possible values of $M$, i.e.,
$0,\dots,\ell-1$. In the first column, we write the possible values of
$\mathsf{a}$, i.e., $0,\dots,\ell-1$  such that $\mathsf{a}$ may take at an edge $e$ of multiplicity $M$. We take $\mathsf{a}(e)=i$ if $i$ satisfies the compatibility condition $\gcd(M,\ell)\mid i$. Then, we fill the slot in the $i$th row and $j$th column of the table with
the corresponding values of $\mathsf{a} \odot M$ if and only if $\mathsf{a}=i$ is compatible with $M=j$.  
For clarity we run the check for $\ell=6$ as well.
$$
\begin{tabular}{|c||c|c|r|l|c|r|}
\hline $\ell=6$& 0& 1 & 2 & 3 & 4 & 5  \\ \hline
       0&   0& 0&  0&  0& 0& 0\\
\hline
       1 &  & 1 & & & & 5\\
\hline
       2 &  & 2& 2& & 4& 4\\
\hline
       3 &  & 3& & 3& & 3\\
\hline
       4 &  & 4& 4& & 2& 2\\
\hline
       5 &  & 5& & & & 1\\
\hline
\end{tabular}
$$
$$
\begin{tabular}{|c||r|r|r|r|r|r|l|l|c|c|r|r|}
\hline $\ell=12$ & 0 & 1 & 2 & 3 & 4 & 5 & 6  & 7 & 8 & 9 & 10 & 11 \\ \hline
       0 & 0& 0& 0& 0& 0& 0& 0& 0& 0& 0& 0& 0 \\
\hline
      1 &  &  1&  & &  &  5&  &  7&  &  &  & 11 \\
\hline
      2 &  &  2&  2& &  &  10&  &  2&  &  &  10& 10 \\
\hline
      3 &  &  3&  & 3&  &  3&  &  9&  &  9&  & 9 \\
\hline
      4 &  &  4&  4& &  4&  8&  &  4&  8&  &  8& 8 \\
\hline
      5 &  &  5&  & &  &  1&  &  11&  &  &  & 7 \\
\hline
      6 &  &  6&  6& 6&  &  6&  6&  6&  &  6&  6& 6 \\
\hline
      7 &  &  7&  & &  &  11&  &  1&  &  &  & 5 \\
\hline
      8 &  &  8&  8& &  8&  4&  &  8&  4&  & 4& 4 \\
\hline
      9 &  &  9&  & 9&  &  9&  &  3&  & 3&  & 3 \\
\hline
      10 &  &  10&  10& &  &  2&  &  10&  &  &  2& 2 \\
\hline
      11 &  &  11&  & &  &  7&  &  5&  &  &  & 1 \\ \hline
\end{tabular}
$$
\begin{center}
Table 2. Multiplication tables for $\odot$ and $\ell$= 6 and 12.
\end{center}

Notice that we fill the table with the corresponding values of $\mathsf{a}\odot M$, if $\mathsf{a}$ is a multiple of $M$. Without loss of generality, we arrange $a_{1}, a_{2}$ and $a_{3}$ in unordered 3-tuple $(a_{1},a_{2},a_{3})$ and  consider all possible cases such that $a_{1}+a_{2}+a_{3}<\ell=12$. Then, the following cases are obtained:
\begin{center}
(1,2,3), (1,2,4), (1,2,5), (1,2,6), (1,2,7), (1,2,8), (1,3,4), (1,3,5), (1,3,6), (1,3,7), (1,4,5), (1,4,6), (2,3,4), (2,3,5), (2,3,6), (1,1,1), (1,2,1), (1,3,1), (1,4,1), (1,5,1), (1,6,1), (1,7,1), (1,8,1), (1,9,1), (2,1,2), (2,2,2), (2,3,2), (2,4,2), (2,5,2), (2,6,2), (2,7,2), (3,1,3), (3,2,3), (3,3,3), (3,4,3), (3,5,3), (4,1,4), (4,2,4), (4,3,4), (5,1,5).
\end{center}
 On the other hand the action of automorphism $\mathsf{a}$ on $M$ should be in the image of $\delta$. Therefore,  we look for  two numbers in the table such that they are equal and we also look for the two numbers in the table such that the sum of them is $\ell=12$ simultaneously. Then,  the following cases of automorphism $\mathsf{a}$ remain:
\begin{center}
(1,1,1), (1,5,1), (1,7,1), (2,2,2), (3,3,3), (5,1,5).
\end{center}
If take $(1,1,1)$, then we can choose $m_{1}, m_{2}$ and $m_{3}$ as the following list:\\
$(m_{1}=1,m_{2}=11,m_{3}=1)$,  $(m_{1}=11,m_{2}=1,m_{3}=11)$, $(m_{1}=7,m_{2}=5,m_{3}=7)$ and $(m_{1}=5,m_{2}=7,m_{3}=5)$. \\
If take $(1,5,1)$ or $(3,3,3)$ or $(5,1,5)$, then we can choose $m_{1}, m_{2}$ and $m_{3}$ as the following list:\\
$(m_{1}=1,m_{2}=7,m_{3}=1)$,  $(m_{1}=7,m_{2}=1,m_{3}=7)$, $(m_{1}=5,m_{2}=11,m_{3}=5)$ and $(m_{1}=11,m_{2}=5,m_{3}=11)$. \\
 Finally If take $(1,7,1)$ or $(2,2,2)$, then we can choose $m_{1}, m_{2}$ and $m_{3}$ as the following list:\\
$(m_{1}=1,m_{2}=5,m_{3}=1), (m_{1}=5,m_{2}=1,m_{3}=5), (m_{1}=7,m_{2}=11,m_{3}=7)$ and  $(m_{1}=11,m_{2}=7,m_{3}=11)$.\\
For these remaining cases, $m_{1}+m_{3}\neq m_{2}$. It means that $M$ doesn't lie in  the $\mathrm{Ker} \ \partial$. This implies the claim that  there is no junior ghost $\mathsf{a}$  for any stable graph for which $M \in \mathrm{Ker} \ \partial$ and $\mathsf{a}\odot M \in  \mathrm{Im} \ \delta$ for $\ell=12$.

Now we want to show there exists a codimension 4 locus of noncanonical singularities for $\ell=12$. Hence it suffices to take  $m_{1}= 1, m_{2}=5, m_{3}=m_{4}=2$ and  $a_{1}=a_{2}=a_{3}=a_{4}=2$. Then, there exists a junior ghost automorphism.
$$
\begin{matrix}
\xymatrix{
\bullet \ar@/^2pc/[rr]^{m_1=1} \ar@/
_1.5pc/[rr]^{m_3=2} \ar@/_2pc/[rr]_{m_4=2}
& & \bullet \ar[ll]_{m_2=5}
}
\end{matrix}
$$
\begin{center}
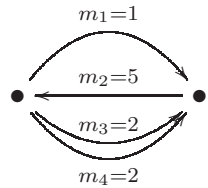
Figure 3. Graph $M$ with two vertices and four edges.
\end{center}

\textbf{Acknowledgements.} This work has been initiated in Pragmatic 2015. I am grateful to Alfio Ragusa, Francesco Russo and Giuseppe Zappal\`a for organising this summer school. The paper focuses on an open problem discussed in the paper titled "Singularities of the Moduli Space of Level Curves"  written by Alessandro Chiodo and Gavril Farkas.  I am extremely grateful to Alessandro Chiodo for useful explanations and support. I also appreciate Hassan Haghighi for his advice.


(S. Tashvighi) \textsc{Faculty  of Mathematics K.N. Toosi University
of Technology, Tehran, Iran}
\\
\textit{E-mail address}: \texttt{sepid.tashvighi@email.kntu.ac.ir}
\\

\end{document}